\documentstyle[12pt,amsfonts]{article}

\input{epsf.tex}

\newlength{\grafikwidth}
\newlength{\twografikwidth}
\newlength{\threegrafikwidth}
\newcommand{\Torus}{\begin{figure}
	\centerline{
		\epsfxsize=\grafikwidth\framebox{\epsfclipon\epsfbox{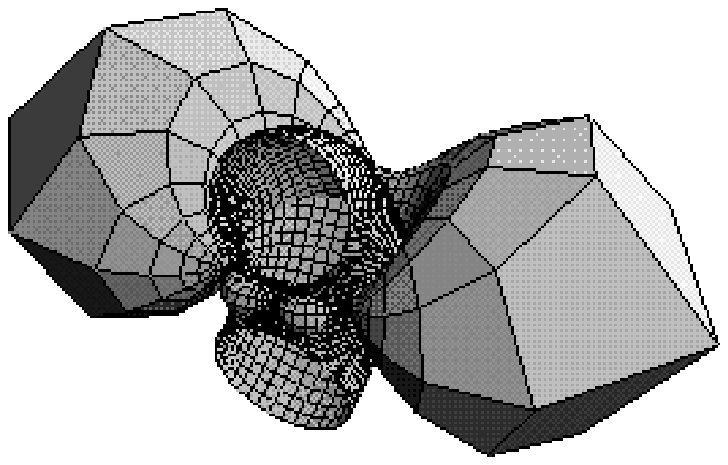}} }
	\caption{A typical Darboux transform of the Clifford Torus}
	\label{FTorus}\end{figure}}
\newcommand{\IvanOne}{\begin{figure}
	\centerline{
		\epsfxsize=\threegrafikwidth\framebox{\epsfclipon\epsfbox{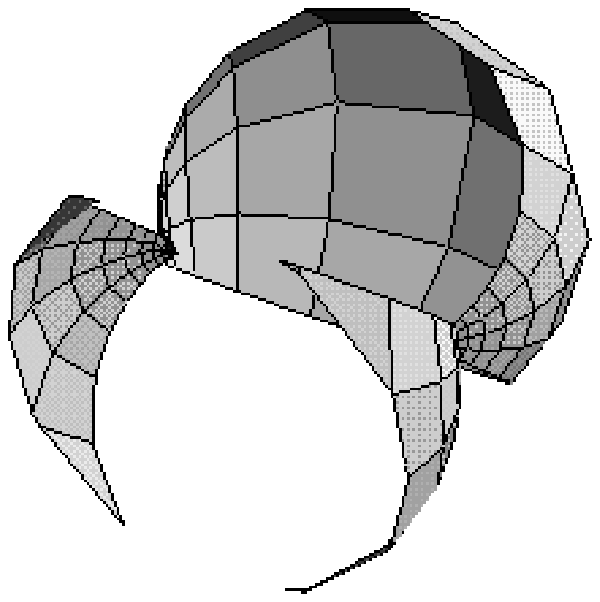}}
		\epsfxsize=\threegrafikwidth\framebox{\epsfclipon\epsfbox{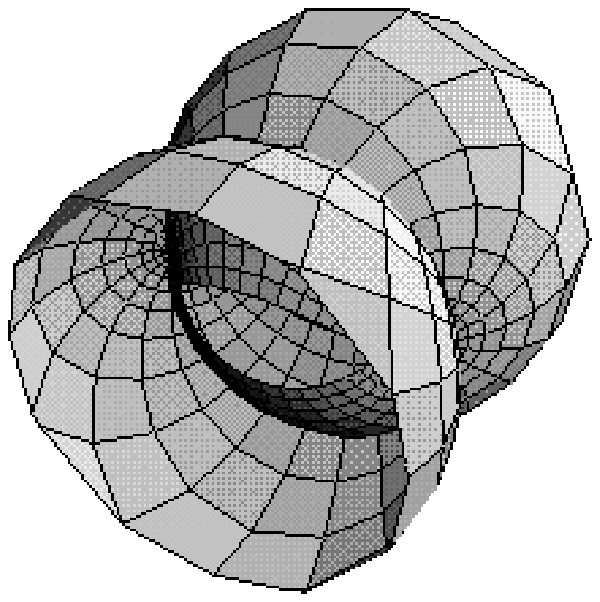}}
		\epsfxsize=\threegrafikwidth\framebox{\epsfclipon\epsfbox{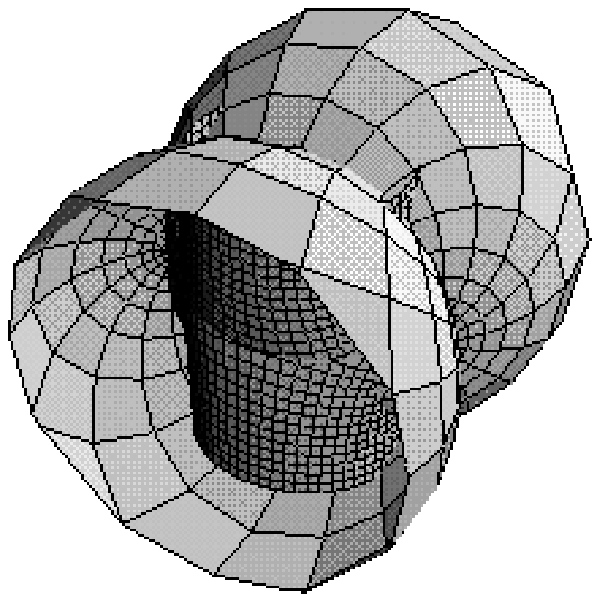}} }
	\vspace{1ex}\centerline{
		\epsfxsize=\grafikwidth\framebox{\epsfclipon\epsfbox{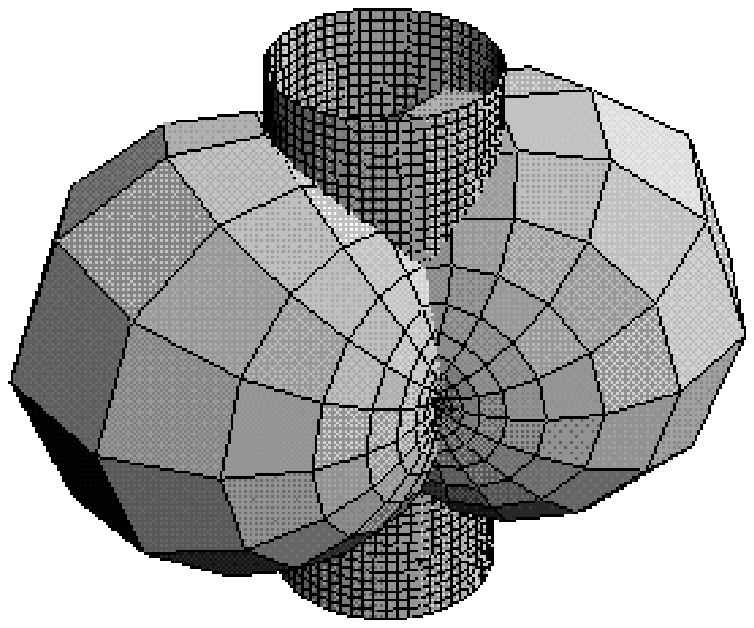}} }
	\caption{A Darboux transform of the cylinder}
	\label{FIvanOne}\end{figure}}
\newcommand{\IvanTwo}{\begin{figure}
	\centerline{
		\epsfxsize=\grafikwidth\framebox{\epsfclipon\epsfbox{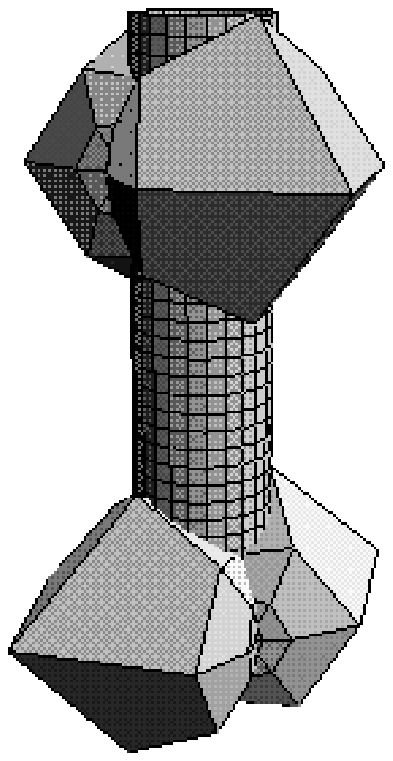}} }
	\caption{A double Darboux transform of the cylinder}
	\label{FIvanTwo}\end{figure}}
\newcommand{\CatMin}{\begin{figure}
	\centerline{
		\epsfxsize=\twografikwidth\framebox{\epsfclipon\epsfbox{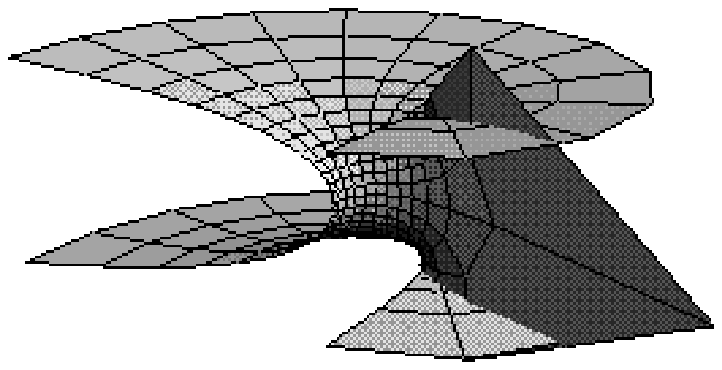}}
		\epsfxsize=\twografikwidth\framebox{\epsfclipon\epsfbox{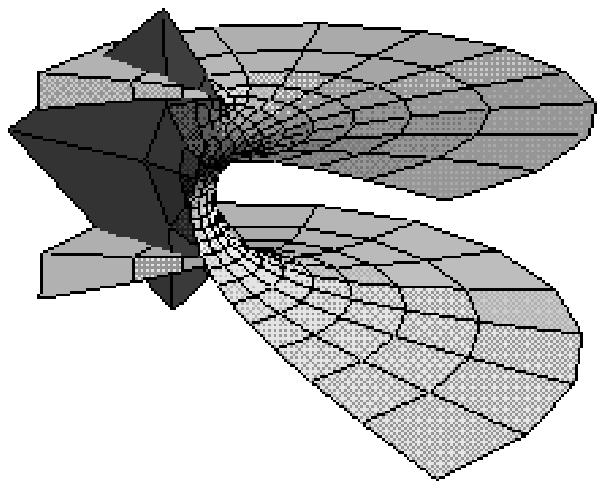}} }
	\vspace{1ex}\centerline{
		\epsfxsize=\grafikwidth\framebox{\epsfclipon\epsfbox{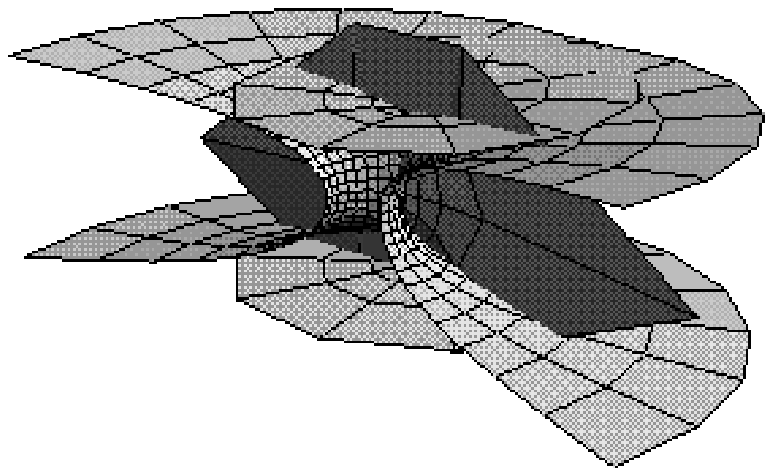}} }
	\caption{A minimal Darboux transform of the Catenoid}
	\label{FCatMin}\end{figure}}
\newcommand{\CatH}{\begin{figure}
	\centerline{
		\epsfxsize=\grafikwidth\framebox{\epsfclipon\epsfbox{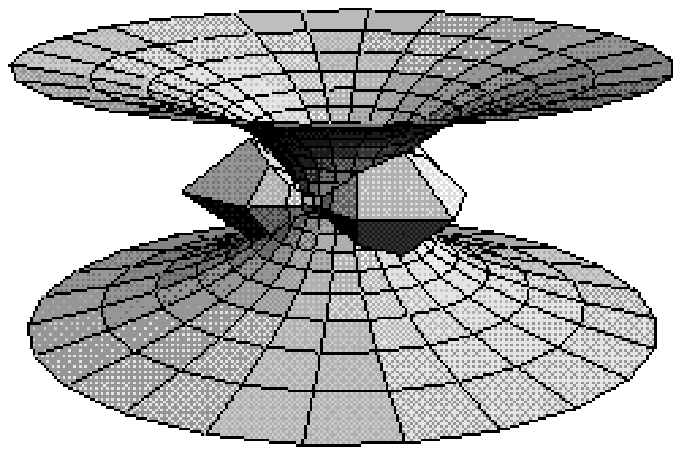}} }
	\caption{A Darboux transform of the Catenoid}
	\label{FCatH}\end{figure}}
\newcommand{\CatV}{\begin{figure}
	\centerline{
		\epsfxsize=\grafikwidth\framebox{\epsfclipon\epsfbox{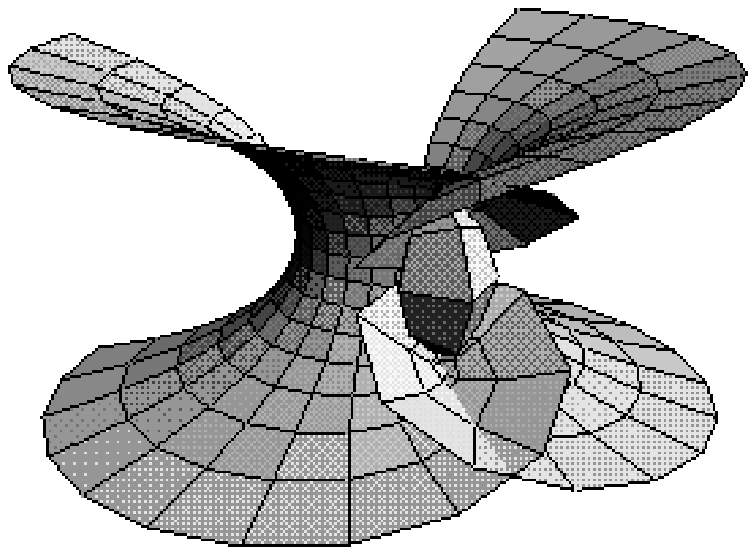}} }
	\caption{Another Darboux transform of the Catenoid}
	\label{FCatV}\end{figure}}

\newcommand{\df}{\sc}

\def\Z{\hspace*{.45ex}\makebox[-.45ex][l]{\hbox{\raise.38ex\hbox{\scriptsize\sl/}}}\rm Z}
\def\Rationals{\hspace*{.26ex}\makebox[-.26ex][l]{\hbox{\raise .53ex\hbox{\tiny/}}}Q}
\def\R{I\hspace*{-.9ex}R}
\def\C{\hspace*{.26ex}\makebox[-.26ex][l]{\hbox{\raise .53ex\hbox{\tiny/}}}C}
\def\H{I\hspace*{-.9ex}H}
\renewcommand{\Im}{{\rm Im}}
\renewcommand{\Re}{{\rm Re}}
\newcommand{\Q}{{\rm Q}}
\newcommand{\DV}{{\rm D\!V}}
\newcommand{\widebar}[1]{\overline{#1}}

\newcommand{\Theorem}[2]{\paragraph{#1:}{\em#2\vspace{2ex}}}

\title{A discrete version of the Darboux transform for isothermic surfaces}
\author{Udo Hertrich-Jeromin\thanks{Partially supported by the Alexander von
	Humboldt Stiftung and by the Forschungsinstitut f\"ur
	Mathematik, ETH Z\"urich}, Tim Hoffmann and Ulrich Pinkall}

\begin{document}
\maketitle

\begin{abstract} We study Christoffel and Darboux transforms of discrete
	isothermic nets in 4-dimensional Euclidean space: definitions
	and basic properties are derived. Analogies with the smooth
	case are discussed and a definition for discrete Ribaucour
	congruences is given.
	Surfaces of constant mean curvature are special among all
	isothermic surfaces: they can be characterized by the fact
	that their parallel constant mean curvature surfaces are
	Christoffel and Darboux transforms at the same time. This
	characterization is used to define discrete nets of constant
	mean curvature. Basic properties of discrete nets of constant
	mean curvature are derived.
	\end{abstract}

\section{Introduction}
Stimulated by the integrable system approach to isothermic surfaces (cf.
\cite{BJPP} or \cite{Cieslinski}) --- or, better: to Darboux {\em pairs}
of isothermic surfaces\footnote{Transformations seem to play an important
role in the relation of surface theory and integrable system theory.} ---
``discrete isothermic surfaces'' or ``nets'' have been introduced by
Alexander Bobenko and Ulrich Pinkall \cite{Sasha}.
For the development of this theory it seems having been crucial {\em not}
to use the standard calculus for M\"obius geometry but a quaternionic approach
--- which is well developed in case of Euclidean ambient space.
In \cite{Sasha} the Christoffel transform (or ``dual'') \cite{Christoffel}
for discrete isothermic nets in Euclidean space --- the Christoffel transform
of an isothermic surface is closely related to the Euclidean geometry of the
ambient space --- is developed.
The Christoffel transform is an important tool in the theory of isothermic
surfaces in Euclidean space:
it can be used to characterize minimal surfaces (cf. \cite{Christoffel},
\cite{Sasha}) and surfaces of constant mean curvature (cf. \cite{Q3}) between
all isothermic surfaces --- for discrete nets of constant mean curvature we
will later use this characterization as a definition.

Recently, a quaternionic calculus for M\"obius differential geometry \cite{Q1}
was developed. This approach turned out to be very well adapted to the theory
of (Darboux pairs of) isothermic surfaces \cite{Q2} --- based on these results
a quaternionic description for the Darboux transform of isothermic surfaces
in Euclidean 4-space $\R^4\cong\H$ was given \cite{Q3}:
any Darboux transform of an isothermic surface in $\R^4$ can be obtained as
the solution of a Riccati type partial differential equation.
This equation can easily be discretized which provides us with a discrete
version of the Darboux transform.
Note that it seems to be more natural to work in 4-dimensional ambient space
rather than in the codimension 1 setting\footnote{Apart from the fact that
calculations can be done much easier in 3- or 4-dimensional ambient space
it seems to be clear that most considerations on the Darboux transform of
discrete isothermic nets can be done in arbitrary dimensions --- ironically,
there might occur problems with the Christoffel transform.}:
for example the structure of our Riccati type equation becomes more clear
in 4-dimensional ambient space.

After a comprehensive discussion of the cross ratio in Euclidean 4-space
we recall the definition of a discrete isothermic net and prove the basic
facts on the Christoffel transform in 4-dimensional ambient space.
Then, we give a definition of the Darboux transform of discrete isothermic
nets by discretizing the Riccati type partial differential equation which
defines the Darboux transforms of an isothermic surface in the smooth case.
In the smooth setting a Darboux transform of an isothermic surface is
obtained as the second envelope of a suitable 2-parameter family of
2-spheres. In the discrete setting there naturally appears a ``2-parameter
family'' of 2-spheres, too --- this suggests a definition of the
``discrete envelopes'' of a ``discrete Ribaucour congruence'' which we
briefly discuss. 
Following the definitions we prove Bianchi's permutability theorems
\cite{Bianchi2} for multiple Darboux transforms and for the Darboux and the
Christoffel transform in the discrete case.
In the final section we discuss discrete nets of constant mean curvature:
(smooth) cmc surfaces can be characterized by the fact that their
parallel cmc surface is a Christoffel and a Darboux transform of
the original surface at the same time.
We use this as a definition and discuss Darboux transforms of constant
mean curvature of discrete nets of constant mean curvature.

To illustrate the effect of the Darboux transformation on discrete isothermic
nets we have included a couple of pictures.
Comparing these pictures with those in \cite{Q3} (which were calculated using
the ``smooth theory'') the reader will observe considerable similarities ---
indicating another time the close relation with the ``smooth theory''.
A discrete isothermic net in Euclidean 3-space allows $\infty^4$ Darboux
transforms into discrete isothermic nets in 3-space. Figure \ref{FTorus}
shows a Darboux transform of a discrete isothermic net on the Clifford
torus which is closed in one direction. The reader will recognize a typical
behaviour of Darboux transforms: along one curvature line ``bubbles'' are
``added'' to the original net while ,in the other direction, the transform
approaches the original net asymptotically.
From this typical behaviour we might expect that Darboux transforms of
(smooth or discrete) isothermic nets on a torus will never become doubly
periodic.

Figure \ref{FIvanOne} shows a periodic constant mean curvature Darboux transform
of a discrete isothermic net on the cylinder --- the correspoding smooth Darboux
transform is a B\"acklund transform of the cylinder, at the same time \cite{Q3}.
The (also periodic) isothermic net (of constant mean curvature) shown in figure
\ref{FIvanTwo} is obtained by applying a second Darboux transformation.
In the last section of the present paper we will show that discrete isothermic
nets of constant mean curvature $H\neq0$ allow $\infty^3$ Darboux transforms of
constant mean curvature $H$.
A similar theorem holds for discrete minimal nets\footnote{Discrete minimal
nets are defined in \cite{Sasha}. As the proof for the theorem on cmc Darboux
transforms of cmc nets relies on the permutability theorem for multiple Darboux
transforms, the proof for that on minimal Darboux transforms of minimal nets
depends on the permutability theorem for Christoffel and Darboux transforms
--- when we characterize minimal isothermic nets by the fact that their
Christoffel transforms are isothermic nets on a 2-sphere.} (cf. \cite{Q3})
--- a minimal Darboux transform of an isothermic net on the catenoid is shown
in figure \ref{FCatMin}.
Other (obviously non minimal) Darboux transforms of the same catenoid net
with ``bubbles added'' along different directions are displayed in figures
\ref{FCatH} and \ref{FCatV} --- these surfaces might help to ``see'' the
``typical behaviour'' of the minimal Darboux transform in figure \ref{FCatMin}
...

\Torus 

Since the cross ratio plays a key role in our investigations on discrete
isothermic surfaces and the Darboux transform it seems useful to discuss
properties of
\section{The Cross Ratio} of four points in Euclidean 4-space.
In \cite{Sasha} Alexander Bobenko and Ulrich Pinkall introduce the (complex)
cross ratio $$\DV(Q_1,Q_2,Q_3,Q_4)=\Re\Q(Q_1,Q_2,Q_3,Q_4)\pm
i\cdot|\Im\Q(Q_1,Q_2,Q_3,Q_4)|$$ where $\Q(Q_1,Q_2,Q_3,Q_4)=
(Q_1-Q_2)(Q_2-Q_3)^{-1}(Q_3-Q_4)(Q_4-Q_1)^{-1}$ for four points
$Q_1,Q_2,Q_3,Q_4\in\R^3\cong\Im\H$ in Euclidean 3-space.
If $S\subset\R^3$ is a 2-sphere containing these four points this is exactly
the cross ratio of the four points interpreted as complex numbers on $S$ as the
Riemann sphere\footnote{Note that this cross ratio is slightly different from
the classical cross ratio: the classical one may be obtained by interchanging
$Q_2$ and $Q_3$ in the above formula. This reflects the fact that in discrete
surface theory it is more useful to consider the cross ratio as an invariant
of a surface patch rather than an invariant of two point pairs.}
--- the ambiguity of the sign of the imaginary part in the cross ratio
corresponds to the possible orientations of the sphere.
Note that the imaginary part of the cross ratio vanishes exactly when the four
points are concircular\footnote{Or, they are colinear. --- But since we are doing
M\"obius geometry we do not distinguish circles and straight lines.}. In
this case the four points do not determine a unique 2-sphere but a whole
pencil of 2-spheres.

This construction holds perfectly for four points $Q_i$, $i=1\dots4$, in
4-dimensional Euclidean space $\R^4\cong\H$:

\Theorem{Definition (Cross Ratio)}{The complex number
\begin{equation}\label{CrossRatio}\hspace{-1.29em}\begin{array}{rl}
	\DV(Q_1,Q_2,Q_3,Q_4):=&
		\Re\Q(Q_1,Q_2,Q_3,Q_4)+i\cdot|\Im\Q(Q_1,Q_2,Q_3,Q_4)|,\cr
	\Q(Q_1,Q_2,Q_3,Q_4):=&
		(Q_1-Q_2)(Q_2-Q_3)^{-1}(Q_3-Q_4)(Q_4-Q_1)^{-1}\cr
	\end{array}\end{equation}
is called the {\df cross ratio} of four points $Q_1,Q_2,Q_3,Q_4\in\R^4\cong\H$
in Euclidean 4-space\footnote{Note that for four points in the complex plane
$\C\subset\H$ (which already carries an orientation) we may simply use
$\Q$ itself as the cross ratio.\label{Orientation}} --- or of the quadrilateral
$(Q_1,Q_2,Q_3,Q_4)\subset\R^4$, respectively.}

It is a straightforward calculation to express the cross ratio in terms
of the distances of the four points\footnote{At this point we would like
to thank Matthias Z\"urcher for helpful discussions.}:
\begin{equation}\label{CrossRatioConformal}
	\DV(Q_1,Q_2,Q_3,Q_4)=\frac{[l_{12}^2l_{34}^2+l_{14}^2l_{23}^2-l_{13}^2l_{24}^2]
		+i\sqrt{-\det((l_{ij}^2)_{i,j=1,\dots,4})}}{2l_{14}^2l_{23}^2}
	\end{equation}
where $l_{ij}:=|Q_i-Q_j|$.
Note that the vertices of a planar quadrilateral are concircular if and only
if the product of the diagonal lengthes equals the sum of the products of the
lengthes of opposite sides\footnote{This is reflected by the fact that
$$-\det((l_{ij}^2)_{i,j=1,\dots,4})=(a+b+c)(-a+b+c)(a-b+c)(a+b-c)$$
where $a=l_{12}l_{34}$, $b=l_{13}l_{24}$ and $c=l_{14}l_{23}$.
Note that all factors are non negative which can be interpreted as triangle
inequalities for the triangle with lengthes $a$, $b$ and $c$.
With this interpretation, the determinant in the cross ratio is 16 times
the squared area of this triangle (Heron's formula).}.
Consequently, the cross ratio of concircular points can be calculated from
the lengthes of two pairs of opposite edges.

%
%

If we consider $\R^4\cong\H=\{(Q,1)\,|\,Q\in\H\}$ as (a subset of) the conformal
4-sphere $\H{}P^1=\{V\cdot\H\,|\,V\in\H^2\}$ the squared norms in this equation
are replaced by a biquadratic function $|Q_i-Q_j|^2={\cal D}((Q_i,1),(Q_j,1))$,
i.e. ${\cal D}(V\lambda,W)={\cal D}(V,W\lambda)=|\lambda|^2{\cal D}(V,W)$ for
$V,W\in\H^2$ and $\lambda\in\H$. Thus the cross ratio of four points is
independent of the choice of homogeneous coordinates of the points and
consequently it is a {\em conformal} invariant \cite{Q1}.

As another consequence of the above equation (\ref{CrossRatioConformal}) we
derive the following identities for the cross ratio --- note that the determinant
in (\ref{CrossRatioConformal}) is invariant under permutations of the four
points:
\begin{equation}\label{CrossRatioId1}\begin{array}{rcl}
	\DV(Q_4,Q_1,Q_2,Q_3)&=&\frac{1}{\ \overline{\DV(Q_1,Q_2,Q_3,Q_4)}\ },\cr
	\DV(Q_1,Q_3,Q_2,Q_4)&=&1-\overline{\DV(Q_1,Q_2,Q_3,Q_4)}.\cr
	\end{array}\end{equation}
The complex conjugation in these equations arises since we defined
the cross ratio to have always a non negative imaginary part.
With the help of these two identities we easily obtain a second set of
identities,
\begin{equation}\label{CrossRatioId2}\begin{array}{rcl}
	\DV(Q_1,Q_2,Q_3,Q_4)&=&\DV(Q_3,Q_4,Q_1,Q_2)\cr
	&=&\DV(Q_2,Q_1,Q_4,Q_3)\cr
	&=&\DV(Q_4,Q_3,Q_2,Q_1).\cr
	\end{array}\end{equation}
The complete set of the twenty four identities for the cross ratio can be
derived using these two sets of identites --- they are symbolized in
figure \ref{CrossRatioIdentities}.
\begin{figure}\centerline{\input{cross.tex}}
	\caption{Identities for the complex cross ratio}
	\label{CrossRatioIdentities}\end{figure}

As a preparation we first give a planar version of our hexahedron lemma
--- which will take a central role throughout this paper:

\Theorem{Lemma}{Given a quadrilateral $(x_1,x_2,x_3,x_4)$ in the complex
plane and a complex number $\lambda\in\C$ there is a unique quadrilateral
$(z_1,z_2,z_3,z_4)$ to each initial point $z_1\in\C$ such that the cross
ratios\footnote{Here, we use $\Q$ as the cross ratio of four points
in the complex plane (cf. footnote \ref{Orientation}).} satisfy
\begin{equation}\begin{array}{cccll}
	\Q(z_1,z_2,z_3,z_4)&=&\Q(x_1,x_2,x_3,x_4)&=:&\mu\cr
	\Q(z_1,z_2,x_2,x_1)&=&\Q(z_3,z_4,x_4,x_3)&=&\mu\lambda\cr
	\Q(z_2,z_3,x_3,x_2)&=&\Q(z_4,z_1,x_1,x_4)&=&\lambda.\cr
	\end{array}\end{equation}}

\begin{figure}\centerline{\input{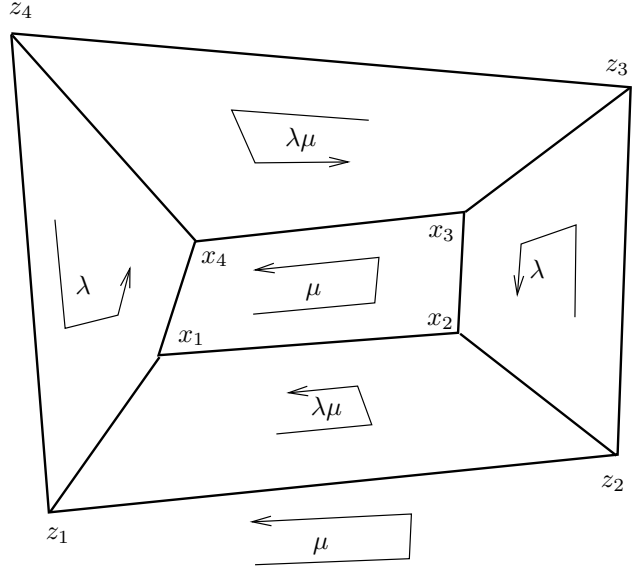}}
	\caption{The hexahedron lemma (planar version)}
	\label{hexahedron}\end{figure}

This lemma may be proven by a straightforward algebraic computation:
first we solve the equations $\Q(z_1,z_2,x_2,x_1)=\mu\lambda$,
$\Q(z_2,z_3,x_3,x_2)=\lambda$ and $\Q(z_3,z_4,x_4,x_3)=\mu\lambda$
succesively for $z_2$, $z_3$ and $z_4$. This defines us a {\em unique}
quadrilateral $(z_1,z_2,z_3,z_4)$ in the complex plane. What is left
to do is to verify the remaining two cross ratios
--- since the calculations become fairly long we prefered to let the computer
algebra system ``Mathematica'' do the work:
\begin{verse}{\tt
	QDVSolve[a$_-$,Q1$_-$,Q2$_-$,Q3$_-$] := Block[\{X\},\\*
	\hspace{1em}X = (Q2-Q3) (Q1-Q2)$^{\wedge}$(-1) a;\\*
	\hspace{1em}(1+X)$^{\wedge}$(-1) (X Q1 + Q3)\\*
	\hspace{1em}];\\*
	(* cross ratio of the initial quadrilateral: *)\\*
	x4 = QDVSolve[ m, x1, x2, x3 ];\\*
	(* we compute the second quadrilateral: *)\\*
	z2 = QDVSolve[ m l, x2, x1, z1 ];\\*
	z3 = QDVSolve[ 1 l, x3, x2, z2 ];\\*
	z4 = QDVSolve[ m l, x4, x3, z3 ];\\*
	(* we verify the remaining two cross ratios: *)\\*
	Factor[ z1 - QDVSolve[ l, x1, x4, z4 ] ]\\*
	Factor[ z4 - QDVSolve[ m, z1, z2, z3 ] ]
	}\end{verse}
This completes the proof of the above lemma --- note, that we could have
omitted the last calculation since the cross ratio of the constructed
quadrilateral can also be derived directly from figure \ref{hexahedron}.

This lemma can not be directly generalized to quadrilaterals in space since
a given (complex) cross ratio does no longer determine the fourth vertex of
a quadrilateral in space. But, in our applications we will only be interested
in {\em real} cross ratios --- and in this case the cross ratio {\em does}
determine the fourth vertex. Thus, if we start with an initial quadrilateral
$(X_1,X_2,X_3,X_4)$ with concircular vertices and one vertex $Z_1$ of the
other vertex we can construct succesively the other vertices $Z_2$, $Z_3$
and $Z_4$ on the circles given by the point triples $\{X_1,X_2,Z_1\}$,
$\{X_2,X_3,Z_2\}$ and $\{X_3,X_4,Z_3\}$. Since these circles are always
given by three points lying on the 2-sphere $S$ containing the circle
through the vertices of the first quadrilateral and the initial vertex
$Z_1$ the whole construction takes place on this 2-sphere. Consequently,
we can prescribe one (real) cross ratio as above:

\Theorem{The Hexahedron lemma}{Given a quadrilateral $(X_1,X_2,X_3,X_4)$ in
Euclidean 4-space with concircular vertices and a real number $\lambda\in\R$
there is a unique quadrilateral $(Z_1,Z_2,Z_3,Z_4)$ to each initial point
$Z_1\in\H$ such that
\begin{equation}\label{HexahedronLemma}\begin{array}{cccll}
	\DV(Z_1,Z_2,Z_3,Z_4)&=&\DV(X_1,X_2,X_3,X_4)&=:&\mu\cr
	\DV(Z_1,Z_2,X_2,X_1)&=&\DV(Z_3,Z_4,X_4,X_3)&=&\mu\lambda\cr
	\DV(Z_2,Z_3,X_3,X_2)&=&\DV(Z_4,Z_1,X_1,X_4)&=&\lambda.\cr
	\end{array}\end{equation}
Moreover, the vertices of the hexahedron $(X_1,X_2,X_3,X_4;Z_1,Z_2,Z_3,Z_4)$
lie on a 2-sphere $S\subset\H$.}

Now, we are prepared to start our investigations on
\section{Discrete isothermic nets}
Away from umbilics isothermic surfaces in Euclidean 3-space can be
defined by the existence of conformal curvature line parameters;
the results from \cite{Q2} and \cite{Q3} suggest that this is a
good definition in case of 4-dimensional ambient space, too.
In the codimension 2 case this definition requires the normal bundle
of the surface to be flat --- otherwise, it makes no sense to speak
of ``curvature line coordinates''.

Following \cite{BoPi} a surface is parametrized by curvature lines
if and only if the parameter lines divide the surface in infinitesimal
(planar) rectangles and the surface is isothermic if and only if it
is divided into infinitesimal squares by its lines of curvature
(cf. \cite{Cayley}).
We may reformulate these criteria in a M\"obius invariant flavour
using the cross ratio:

\Theorem{Lemma}{A (smooth) surface is parametrized by curvature lines
if and only if the parameter lines devide the surface in infinitesimal
patches with negative real cross ratios and\\
the surface is isothermic if and only if its curvature lines devide
the surface in infinitesimal harmonic patches, i.e. with cross ratios -1.}

\CatH 

These criteria may motivate the following definition of discrete isothermic
nets --- note that we prefer the notion of a ``net'' since this term refers
to a {\em parametrized} surface rather than to a surface in space as a purely
geometric object: a ``discrete isothermic net'' is the analog of an isothermic
surface given through conformal curvature line parameters.

\Theorem{Definition (Discrete Curvature line Net, Discrete Isothermic Net)}{
A map $F:\Gamma\to\H$, $\Gamma\subset\Z^2$, is called a {\df discrete curvature
line net} if the vertices of all elementary quadrilaterals are concircular and
the quadrilaterals are embedded:
\begin{equation}\label{DiscreteCNet}
	\DV(F_{m,n},F_{m+1,n},F_{m+1,n+1},F_{m,n+1})<0;\end{equation}
it is called an {\df isothermic net} if all elementary quadrilaterals are
harmonic:
\begin{equation}\label{DiscreteINet}
	\DV(F_{m,n},F_{m+1,n},F_{m+1,n+1},F_{m,n+1})=-1.\end{equation}}

It is possible to define more general ``isothermic nets'' or ``isothermic
surfaces'' as a discrete version of isothermic surfaces given in arbitrary
curvature line coordinates (see \cite{Sasha}).
Even though some formulas have to be modified in this approach
(see for example the definition of the ``dual surface'' in \cite{Sasha})
all the facts on the Christoffel and Darboux transforms of isothermic nets
we are going to develop seem to hold in the more general setting.
However, for this presentation we prefer the simpler and more enlightening
approach we just introduced.

As in the codimension 1 case, discrete isothermic nets have a Christoffel
transform (dual). And, in some sense, they can be {\em
characterized} by the existence of a Christoffel transform:

\Theorem{Theorem and Definition (Christoffel transform)}{If a discrete net
$F:\Gamma\to\H$ is isothermic then
\begin{equation}\label{ChristoffelTransform}\begin{array}{rcl}
	1&=&\lambda(F_{m+1,n}-F_{m,n})(F^c_{m+1,n}-F^c_{m,n})\cr
	-1&=&\lambda(F_{m,n+1}-F_{m,n})(F^c_{m,n+1}-F^c_{m,n}),\cr
	\end{array}\end{equation}
$\lambda\in\R\setminus\!\{0\}$, defines another discrete isothermic net
$F^c:\Gamma\to\H$, a {\df Christoffel transform} (or {\df dual}) of $F$.
Note that a Christoffel transform $F^c$ of an isothermic net $F$ is
determined up to a scaling and its position in space.\\
If, on the other hand, (\ref{ChristoffelTransform}) defines a second
discrete net --- i.e. (\ref{ChristoffelTransform}) is integrable ---
then both nets are isothermic, or parallelogram nets.}

This theorem has an analog for smooth surfaces: given a parametrized surface
$F:U\subset\R^2\to\H$ the integrability condition for the $\H$-valued 1-form
$F_x^{-1}dx-F_y^{-1}dy$ is equivalent to $F_{xy}^{\perp}=0$ and $F$ conformal
--- i.e. $F$ is a conformal curvature line parametrization of an isothermic
surface --- or $F_{xy}=0$ and hence $F(x,y)=F_1(x)+F_2(y)$ is a translation
surface --- which is the smooth analog of a parallelogram net.

During the proof of this theorem let us denote by
$a:=F_{m+1,n}-F_{m,n}$, $b:=F_{m+1,n+1}-F_{m+1,n}$, $c:=F_{m,n+1}-F_{m+1,n+1}$
and $d:=F_{m,n}-F_{m,n+1}$ the edges of an elementary quadrilateral of the net
$F$. Clearly, we have
\begin{equation}\label{IntegrabilityF}
	0=a+b+c+d\end{equation}
--- the closing condition for the net $F$.
The closing condition (``integrability'') for the dual net $F^c$ reads
\begin{equation}\label{IntegrabilityFc}\begin{array}{rcl}
	0&=&\frac{1}{a}-\frac{1}{b}+\frac{1}{c}-\frac{1}{d}\cr
	&=&\frac{\bar{a}}{|a|^2}-\frac{\bar{b}}{|b|^2}
		+\frac{\bar{c}}{|c|^2}-\frac{\bar{d}}{|d|^2}.\cr
	\end{array}\end{equation}
As a first consequence this shows that a necessary condition for the net $F^c$
to close up is that the elementary quadrilaterals of $F$ be planar.
If we {\em assume} the quadrilateral of $F$ to be planar we can rewrite
(\ref{IntegrabilityFc}) as
\begin{equation}\label{XXX}\begin{array}{rcl}
	0&=&a\left(\frac{1}{a}-\frac{1}{b}+\frac{1}{c}-\frac{1}{d}\right)c\cr
	&=&-[1+ab^{-1}cd^{-1}]d-[1+ad^{-1}cb^{-1}]b\cr
	&=&[1+ab^{-1}cd^{-1}]\cdot[a+c]\cr
	\end{array}\end{equation}
where we used $b\bar{c}d=d\bar{c}b$ for the last equality: see (\ref{Determinant}).

If $F$ is an isothermic net its elementary quadrilaterals are planar with
cross ratio $ab^{-1}cd^{-1}=-1$. Consequently, the net $F^c$ closes up.
The only thing left is to show that elementary quadrilaterals of $F^c$
have cross ratio -1. But,
\begin{equation}\label{XXY}\begin{array}{rcccccl}
	-1&=&ab^{-1}cd^{-1}&=&d^{-1}cb^{-1}a&=&[\frac{1}{a}\,b\,\frac{1}{c}\,d]^{-1}\cr
	\end{array}\end{equation}
which shows that $F^c$ is an isothermic net, too. This proves the first part
of our theorem.

To attack the second part we note that from the above calculation
(\ref{XXX}) it follows that two cases can occur if (\ref{ChristoffelTransform})
is integrable: the quadrilateral under investigation has cross ratio
$-1=ab^{-1}cd^{-1}$ or it is a parallelogram.
Assuming that the quadrilaterals of $F$ do not change their type we conclude
that $F$ and (from the first part) $F^c$ are isothermic nets or both nets are
parallelogram nets. This completes the proof of the above theorem
--- provided we happen to prove the following

\Theorem{Lemma}{Two vectors $a,b\in\R^4\cong\H$ are linearly dependent if and
only if \begin{equation}a\bar{b}=b\bar{a};\end{equation}
three vectors $a,b,c\in\H$ are linearly dependent if and only if
\begin{equation}\label{Determinant}a\bar{b}c=c\bar{b}a;\end{equation}
four vectors $a,b,c,d\in\H$ are linearly dependent if and only if
\begin{equation}a\bar{b}c\bar{d}-\bar{d}c\bar{b}a=
	d\bar{c}b\bar{a}-\bar{a}b\bar{c}d.\end{equation}}

The first statement is immediately clear since this simply means that
$a\,\frac{1}{b}\in\R$. To prove the second statement we calculate
$$a\bar{b}c-c\bar{b}a=2[\det(a_1,b_1,c_1),-(a_0\cdot b_1\times c_1+b_0
\cdot c_1\times a_1+c_0\cdot a_1\times b_1)]$$ where $x_0=\Re x$ and
$x_1=\Im x$. A careful analysis of this equation provides the second
statement. And third we have $$\Re[a\bar{b}c\bar{d}-\bar{d}c\bar{b}a]
=2\det(a,b,c,d)$$ which completes the proof\footnote{At this point we
would like to thank Ekkehard Tjaden for helpful discussions.}.

To motivate our ansatz for the Darboux transform of discrete isothermic
nets let us shortly recall some facts on
\section{The Darboux transform}\label{SmoothDarboux}
of smooth isothermic surfaces. A sphere congruence (a 2-parameter family of
spheres) $S:M^2\to\{\mbox{spheres and planes in $\R^3$}\}$ is said to be
``enveloped'' by a surface $F:M^2\to\R^3$ if at each point the surface has first
order contact to the corresponding sphere\footnote{For simplicity reasons we
do not give {\em exact definitions} --- we just want to give the {\em ideas}
of the used terms. For a comprehensive discussion the reader is refered to
the classical book of W.~Blaschke \cite{Blaschke}.}: $F(p)\in S(p)$ and
$d_pF(T_pM)=T_{F(p)}S$. If a sphere congruence $S$ has two envelopes $F$ and
$\hat{F}$ --- which is, generically, the case --- then it establishes a point
to point correspondance between its two envelopes. From the works of Darboux
\cite{Darboux} and Blaschke \cite{Blaschke} we know that if this correspondance
preserves curvature lines --- the sphere congruence is a ``Ribaucour sphere
congruence'' --- and it is conformal, too, then generically\footnote{If the sphere
congruence is, in a certain sense, ``full'', i.e. if it is not a congruence of
planes in some space of constant curvature.} both envelopes are isothermic
surfaces.
In this case the two envelopes are said to form a ``Darboux pair'', one
surface is said to be a ``Darboux transform'' of the other.
These definitions can be generalized to codimension 2 surfaces in 4-dimensional
space by introducing congruences of 2-spheres --- which are not hyperspheres
any more \cite{Q2}.

In \cite{Q3} we derived a Riccati type partial differential equation
\begin{equation}\label{SmoothRiccati}
	d\hat{F}=(\hat{F}-F)dF^c(\hat{F}-F)\end{equation}
for Darboux transforms $\hat{F}$ of an isothermic surface $F$ where $F^c$
is a Christoffel transform of $F$. Note that in order to obtain {\em every}
Darboux transform $\hat{F}$ of $F$ as a solution of (\ref{Riccati}), it is
crucial not to fix the scaling of the Christoffel transform: rescalings of
$F^c$ make $\hat{F}$ run through the associated family of Darboux transforms.
Considering this Riccati equation as an initial value problem we see that an
isothermic surface allows $\infty^5$ Darboux transforms --- or, if we are
interested in surfaces in 3-dimensional space, it allows $\infty^4$ Darboux
transforms.

\CatV 

This Riccati type equation (\ref{SmoothRiccati}) can now be easily discretized
to obtain a system of difference equations
$$\begin{array}{rcl}
	\lambda(\hat{F}_{m+1,n}-\hat{F}_{m,n})&=&
	(\hat{F}_{m,n}-F_{m,n})(F_{m+1,n}-F_{m,n})^{-1}(\hat{F}_{m+1,n}-F_{m+1,n})\cr
	\lambda(\hat{F}_{m,n+1}-\hat{F}_{m,n})&=&
	-(\hat{F}_{m,n}-F_{m,n})(F_{m,n+1}-F_{m,n})^{-1}(\hat{F}_{m,n+1}-F_{m,n+1})\cr
	\end{array}$$
where we replaced the Christoffel transform $F^c$ of $F$ according to
(\ref{ChristoffelTransform}). However, this ansatz is not unique: we also
could have interchanged the roles of the edges connecting $F$ and $\hat{F}$
by replacing the previous equations by
$$\begin{array}{rcl}
	\lambda(\hat{F}_{m+1,n}-\hat{F}_{m,n})&=&
	(\hat{F}_{m+1,n}-F_{m+1,n})(F_{m+1,n}-F_{m,n})^{-1}(\hat{F}_{m,n}-F_{m,n})\cr
	\lambda(\hat{F}_{m,n+1}-\hat{F}_{m,n})&=&
	-(\hat{F}_{m,n+1}-F_{m,n+1})(F_{m,n+1}-F_{m,n})^{-1}(\hat{F}_{m,n}-F_{m,n})\cr
	\end{array}$$
--- or we even could have used a mean value of the two. But, rewriting these
equations with the cross ratio we see from (\ref{CrossRatioId2}) that both
ansatzes are equivalent.
Thus, without loss of generality we may use our first ansatz to obtain the
following

\Theorem{Theorem and Definition (Darboux Transform)}{If $F:\Gamma\to\H$ is an
isothermic net then the Riccati type system
\begin{equation}\label{Riccati}\begin{array}{rcl}
	\lambda&=&\DV(F_{m,n},\hat{F}_{m,n},\hat{F}_{m+1,n},F_{m+1,n})\cr
	-\lambda&=&\DV(F_{m,n},\hat{F}_{m,n},\hat{F}_{m,n+1},F_{m,n+1})\cr
	\end{array}\end{equation}
defines an isothermic net $\hat{F}:\Gamma\to\H$. Any solution $\hat{F}$ of
(\ref{Riccati}) is called a {\df Darboux transform} of $F$.}

This theorem is an easy consequence of the hexahedron lemma (page
\pageref{HexahedronLemma}). Note that from the identities (\ref{CrossRatioId2})
for the cross ratio we see that the equations (\ref{Riccati}) are symmetric in
$F$ and $\hat{F}$ --- consequently, we also may refer to the pair $(F,\hat{F})$
of isothermic nets as a {\df Darboux pair}.

As a second consequence from the hexahedron lemma we conclude that for a
Darboux pair $(F,\hat{F})$ the vertices of the hexahedron $$(F_{m,n},F_{m+1,n},
F_{m+1,n+1},F_{m,n+1};\hat{F}_{m,n},\hat{F}_{m+1,n},\hat{F}_{m+1,n+1},
\hat{F}_{m,n+1})$$ lie on a sphere $S_{(m,n)^{\ast}}$ --- these spheres
do naturally live on the ``dual lattice''
\begin{equation}\label{DualLattice}\begin{array}{rcl}
	\Gamma^{\ast}&:=&\{((m,n),(m+1,n),(m+1,n+1),(m,n+1))\,|\cr
		&&(m,n),(m+1,n),(m+1,n+1),(m,n+1)\in\Gamma\}\cr
	\end{array}\end{equation}
of $\Gamma$ which consists of the elementary quadrilaterals. Obviously, an
(interior) point pair $(F_{m,n},\hat{F}_{m,n})$ can be obtained as the
intersection of four spheres $S_{(m,n)^{\ast}}$, $S_{(m-1,n)^{\ast}}$,
$S_{(m-1,n-1)^{\ast}}$ and $S_{(m,n-1)^{\ast}}$.
These facts suggest that the Darboux pair $(F,\hat{F})$ ``envelopes''
a discrete Ribaucour congruence:

\Theorem{Definition (Envelopes of a Discrete Ribaucour Congruence)}{
Two discrete curvature line nets $F,\hat{F}:\Gamma\to\R^4$ are said to
{\df envelope a discrete Ribaucour congruence}
$S:\Gamma^{\ast}\to\{\mbox{2-spheres and 2-planes in $\R^4$}\}$
if (interior) point pairs $(F_{m,n},\hat{F}_{m,n})$ lie on four consecutive
spheres $S_{(m,n)^{\ast}}$, $S_{(m-1,n)^{\ast}}$, $S_{(m-1,n-1)^{\ast}}$
and $S_{(m,n-1)^{\ast}}$.}

Note that we require the two envelopes to be curvature line nets
--- and consequently the sphere congruence to be Ribaucour.
If we would omit this assumption the patches of {\em one} envelope could
determine the spheres of the enveloped congruence uniquely since the four
vertices of an elementary quadrilateral would not be concircular in general.
But this would contradict the fact that in the smooth case there is always
a pointwise 1-parameter freedom for an enveloped sphere congruence when one
envelope is given.

The previous discussions indicate a strong similarity to the smooth case
--- to complete this picture we prove two fundamental ``permutability theorems''
which hold in the smooth case:

\Theorem{Theorem}{If $\hat{F}_1,\hat{F}_2:\Gamma\to\H$ are two Darboux
transforms of an isothermic net $F:\Gamma\to\H$ with parameters $\lambda_1$
and $\lambda_2$, respectively, then there exists an isothermic net
$\hat{F}=\hat{F}_{12}=\hat{F}_{21}:\Gamma\to\H$ which is a $\lambda_2$-Darboux
transform of $\hat{F}_1$ and a $\lambda_1$-Darboux transform of $\hat{F}_2$
at the same time.
The nets $F$, $\hat{F}_1$, $\hat{F}_2$ and $\hat{F}$ have constant cross ratio
\begin{equation}\label{CrossRatioPerm}\begin{array}{c}
	\frac{\lambda_2}{\lambda_1}\equiv\DV(F,\hat{F}_2,\hat{F},\hat{F}_1).
	\end{array}\end{equation}}

\begin{figure}\centerline{\input{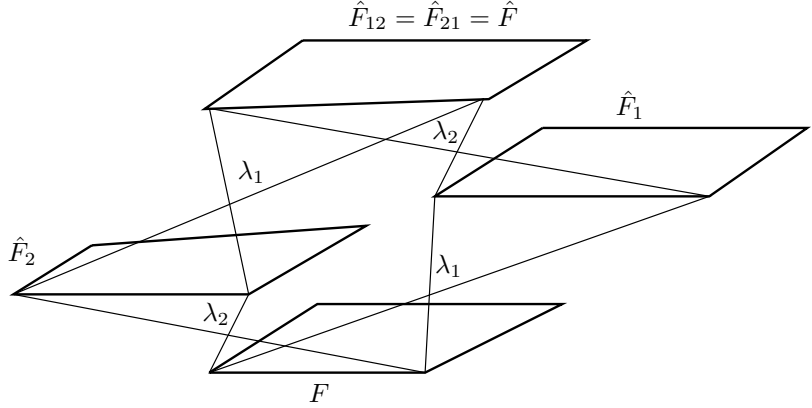}}
	\caption{Bianchi's permutability theorem}
	\label{BianchiTheorem}\end{figure}
This theorem is a discrete analog of Bianchi's permutability theorem
\cite{Bianchi2}. Again, it can be obtained as a consequence of the
hexahedron lemma:
let us pick one edge of an elementery quadrilateral of $F$ and the
corresponding edges of $\hat{F}_1$ and $\hat{F}_2$, as indicated in
figure \ref{BianchiTheorem}.
We denote the vertices of these edges by $X,Y$, $\hat{X}_1,\hat{Y}_1$
and $\hat{X}_2,\hat{Y}_2$. Since $\hat{F}_1$ and $\hat{F}_2$ are Darboux
transforms of $F$ we have
$$\begin{array}{rcl}\DV(X,\hat{X}_1,\hat{Y}_1,Y)&=&\pm\lambda_1,\cr
	\DV(X,\hat{X}_2,\hat{Y}_2,Y)&=&\pm\lambda_2.\cr\end{array}$$
Then, we can find two points $\hat{X},\hat{Y}$ such that
$$\begin{array}{rcccl}\DV(X,\hat{X}_2,\hat{X},\hat{X}_1)&=&
	\DV(Y,\hat{Y}_2,\hat{Y},\hat{Y}_1)&=&\frac{\lambda_2}{\lambda_1}\cr
	\end{array}$$
and, by the hexahedron lemma,
$$\begin{array}{rcl}\DV(\hat{X}_1,\hat{X},\hat{Y},\hat{Y}_1)&=&\pm\lambda_2,\cr
	\DV(\hat{X}_2,\hat{X},\hat{Y},\hat{Y}_2)&=&\pm\lambda_1.\cr\end{array}$$
Repeating this construction for all four edges --- note that by construction the
resulting edges close up to form a quadrilateral --- we obtain an elementary
quadrilateral of an isothermic net $\hat{F}$ which is a $\lambda_1$-Darboux
transform of $\hat{F}_2$ and a $\lambda_2$-Darboux transform of $\hat{F}_1$,
as desired. This completes the proof of our first permutability theorem.

Since the four surfaces in the permutability theorem have constant cross
ratio (\ref{CrossRatioPerm}) it becomes clear that there is an even fancier
version of the permutability theorem which involves not just four but eight
surfaces (cf. \cite{Bianchi2}):
if we start with three Darboux transforms $\hat{F}_1$, $\hat{F}_2$ and
$\hat{F}_3$ of an isothermic net $F$ we can construct three isothermic
nets $\hat{F}_{12}$, $\hat{F}_{23}$ and $\hat{F}_{31}$ via
\begin{equation}\begin{array}{rcl}
	\frac{\lambda_2}{\lambda_1}&\equiv&
		\DV(F,\hat{F}_2,\hat{F}_{12},\hat{F}_1),\cr
	\frac{\lambda_3}{\lambda_2}&\equiv&
		\DV(F,\hat{F}_3,\hat{F}_{23},\hat{F}_2),\cr
	\frac{\lambda_1}{\lambda_3}&\equiv&
		\DV(F,\hat{F}_1,\hat{F}_{31},\hat{F}_3).\cr
	\end{array}\end{equation}
Now, we may apply the construction a second time to obtain an eighth
net $\hat{F}$ which satisfies
\begin{equation}\begin{array}{rcl}    
	\frac{\lambda_2}{\lambda_1}&\equiv&
		\DV(\hat{F}_3,\hat{F}_{31},\hat{F},\hat{F}_{23}),\cr 
	\frac{\lambda_3}{\lambda_2}&\equiv&
		\DV(\hat{F}_1,\hat{F}_{12},\hat{F},\hat{F}_{31}),\cr 
	\frac{\lambda_1}{\lambda_3}&\equiv&
		\DV(\hat{F}_2,\hat{F}_{23},\hat{F},\hat{F}_{12}).\cr 
	\end{array}\end{equation} 
It is a consequence of the hexahedron lemma that the net $\hat{F}$ exists.
This seems to be worth formulating a

\Theorem{Corollary}{If $\hat{F}_i:\Gamma\to\H$, $i=1,2,3$, are Darboux
transforms of an isothermic net $F:\Gamma\to\H$ with parameters $\lambda_i$,
then there exist four isothermic nets $\hat{F}_{ij},\hat{F}:\Gamma\to\H$,
$ij=12,23,31$, such that each $\hat{F}_{ij}$ is a $\lambda_i$-Darboux transform
of $\hat{F}_j$ and a $\lambda_j$-Darboux transform of $\hat{F}_i$ and such that
$\hat{F}$ is a $\lambda_1$-Darboux transform of $\hat{F}_{23}$, a
$\lambda_2$-Darboux transform of $\hat{F}_{31}$ and a $\lambda_3$-Darboux
transform of $\hat{F}_{12}$, at the same time.}

The second ``permutability theorem'' is a discrete version of the compatibility
theorem for the Darboux and Christoffel transform of an isothermic surface
\cite{Q3} (cf. \cite{Bianchi2}):

\Theorem{Theorem}{If $F,\hat{F}:\Gamma\to\H$ form a Darboux pair, then their
Christoffel transforms $F^c,\hat{F}^c=\widehat{F^c}:\Gamma\to\H$ form ---
if correctly scaled and positioned --- a Darboux pair, too.}

To prove this theorem we imitate the proof we gave in the smooth case \cite{Q3}:
first, let us consider the difference functions $G:=\hat{F}-F:\Gamma\to\H$
and $G^c:=\hat{F}^c-F^c:\Gamma\to\H$. Since $\hat{F}$ is a Darboux transform
of $F$ we have
$$\begin{array}{rcl}\lambda(\hat{F}_{m+1,n}-\hat{F}_{m,n})&=&(\hat{F}_{m,n}
	-F_{m,n})(F_{m+1,n}-F_{m,n})^{-1}(\hat{F}_{m+1,n}-F_{m+1,n}),\cr
	\lambda(F_{m+1,n}-F_{m,n})&=&(\hat{F}_{m,n}-F_{m,n})(\hat{F}_{m+1,n}
	-\hat{F}_{m,n})^{-1}(\hat{F}_{m+1,n}-F_{m+1,n}).\cr\end{array}$$
Subtracting the first equation from the second one yields
$$\lambda\left(\frac{1}{\hat{F}_{m+1,n}-F_{m+1,n}}-\frac{1}{\hat{F}_{m,n}
	-F_{m,n}}\right)=\left(\frac{1}{\hat{F}_{m+1,n}-\hat{F}_{m,n}}
	-\frac{1}{F_{m+1,n}-F_{m,n}}\right)$$
and consequently
$$G^c_{m+1,n}-G^c_{m,n}=\frac{1}{G_{m+1,n}}-\frac{1}{G_{m,n}}$$
if $F^c$ and $\hat{F}^c$ are scaled as in (\ref{ChristoffelTransform}).
In a similar way we obtain the same difference equation for the other
direction
$$G^c_{m,n+1}-G^c_{m,n}=\frac{1}{G_{m,n+1}}-\frac{1}{G_{m,n}}.$$
Thus, we can position $F^c$ and $\hat{F}^c$ in $\R^4$ such that
$G^c_{m,n}=G_{m,n}^{-1}.$
As in (\ref{XXY}), we can now calculate the cross ratios
$$\begin{array}{rcl}\DV(F^c_{m,n},\hat{F}^c_{m,n},\hat{F}^c_{m+1,n},
	F^c_{m+1,n})&=&\frac{\lambda^2}{\lambda}=\lambda,\cr
	\DV(F^c_{m,n},\hat{F}^c_{m,n},\hat{F}^c_{m+1,n},F^c_{m+1,n})&
	=&-\lambda\cr\end{array}$$
showing that $F^c$ and $\hat{F}^c$ form indeed a Darboux pair.
This completes the proof.

Troughout the rest of the paper we will consider
\section{Discrete nets of constant mean curvature}
As is well known, (smooth) surfaces of constant mean curvature
(cmc) in 3-dimensional Euclidean space are isothermic surfaces,
i.e. --- away from umbilics --- they allow conformal curvature line
parametrizations.
In fact, we know that surfaces of constant mean curvature ($H\neq0$)
are very special examples of isothermic surfaces in the context of
Darboux and Christoffel transforms:
the (correctly scaled and positioned) Christoffel transform of a
cmc surface $F:U\subset\R^2\to\R^3$ with unit normal field $N:U\to S^2$
is its parallel constant mean curvature surface $F^p=F+\frac{1}{H}N$.
On the other hand, this parallel cmc surface is also a Darboux transform
since a congruence of spheres with constant radius $\frac{1}{2H}$ is
certainly a Ribaucour sphere congruence and --- as is well known ---
a cmc surface $F$ and its parallel cmc surface $F^p$ carry conformally
equivalent metrics.
In \cite{Q3} we proved that this behaviour {\em characterizes} cmc surfaces:
a surface $F:U\to\R^3$ is cmc if and only if its (suitably scaled\footnote{Note
that the sphere --- which is certainly a surface of constant mean curvature ---
behaves exceptionally in this context: any Christoffel transform of the
sphere is a minimal surface but its ``parallel constant mean curvature
surface'' collapses to the center of the sphere. Therefore, for the sphere,
``suitably scaled'' means that the scaling factor $\frac{1}{\lambda}=0$.
\label{Sphere}}
and positioned) Christoffel transform is a Darboux transform of $F$, too.

\CatMin 

Since we are working on Christoffel and Darboux transforms of isothermic
surfaces it seems to be reasonable to use this characterization of cmc
surfaces as a

\Theorem{Definition (Discrete cmc net)}{An isothermic net
$F:\Gamma\to\R^3$ is called a {\df discrete cmc net} if a (suitably
scaled\footnote{Because of the way we introduced the scaling factor
$\lambda$ for the Christoffel transform in (\ref{ChristoffelTransform})
this definition {\em excludes} the sphere: for the sphere the ``suitably
scaled'' Christoffel transform would be the center of the sphere (cf.
footnote \ref{Sphere}), i.e. $\lambda=\infty$.\label{SphereCMC}} and
positioned) Christoffel transform $F^p:\Gamma\to\R^3$ of $F$ is a Darboux
transform of $F$, too.}

However, in Euclidean geometry this definition is rather unsatisfactory
since it doesn't say anything about the mean curvature --- which is claimed
to be constant\footnote{In the geometry of similarities, however, this is
a good definition: here, the value of the (constant) mean curvature is
not an invariant. Note that the geometry of similarities arises much more
naturally (by distinguishing a point at infinity) as a subgeometry of
M\"obius geometry than the Euclidean geometry (where, additionally, a
scaling of the conformal metric has to be choosen).}: recall the meaning
of ``cmc''.
But, according to the above discussions we might define the (constant) mean
curvature\footnote{In our discussions, Alexander Bobenko suggested a definition
for the mean curvature function $H:\Gamma\to\R$ of an (arbitrary) discrete net
$F:\Gamma\to\R^3$ at a point $F_{m,n}$ of the net: consider the point $C$ which
has equal distances $|F_{m-1,n}-C|=|F_{m+1,n}-C|$ and $|F_{m,n-1}-C|=|F_{m,n+1}-C|$
to the pairs of opposite neighbours of $F_{m,n}$ and whose distance to $F_{m,n}$
satisfies $|F_{m,n}-C|^2=\frac{1}{2}(|F_{m+1,n}-C|^2+|F_{m,n+1}-C|^2)$.
Its reciprocal distance $H_{m,n}=|F_{m,n}-C|^{-1}$ can be considered the
mean curvature of the net $F$ at $F_{m,n}$. Note, that $C$ may not be a finite
point --- then the mean curvature vanishes (cf. the definition of discrete minimal
surfaces in \cite{Sasha}).\label{MeanCurv}}
of a discrete cmc net to be the reciprocal of the constant distance (vertex
wise) of the nets $F$ and $F^p$ --- if we happen to prove the following

\Theorem{Theorem}{A Christoffel transform $F^c:\Gamma\to\R^3$ of an
isothermic net $F:\Gamma\to\R^3$ is a Darboux transform of $F$, too, if and
only if the distance $|F^c_{m,n}-F_{m,n}|$ is constant.}

Since the quadrilaterals spanned by corresponding edges\footnote{Corresponding
edges of $F$ and any Christoffel transform $F^c$ are parallel --- note that
we restricted to the codimension 1 case.} of $F$ and $F^c$ are trapezoids this
theorem will be a consequence of

\Theorem{The Trapezium lemma}{The cross ratio of a trapezoid is real
if and only if that trapezoid is isosceles.
A quadrilateral $(Q_1,Q_2,Q_3,Q_4)\subset\H$ is an (isosceles) trapezoid,
i.e. $(Q_1-Q_4)\parallel(Q_2-Q_3)$ (and $|Q_1-Q_2|=|Q_3-Q_4|$), if and only
if\footnote{Note, that if $|Q_1-Q_2|<|Q_1-Q_3|$ then the trapezoid has legs
$Q_1-Q_2$ and $Q_3-Q_4$ and the isosceles trapezoid $(Q_1,Q_2,Q_3,Q_4)$ is
embedded. If $|Q_1-Q_2|>|Q_1-Q_3|$ the edges $Q_1-Q_2$ and $Q_3-Q_4$
become the diagonals of the trapezoid which is then not embedded.}
\begin{equation}\begin{array}{c}
	\DV(Q_1,Q_2,Q_3,Q_4)=-\frac{|Q_1-Q_2|^2}{|Q_1-Q_3|^2-|Q_1-Q_2|^2}\ .
	\end{array}\end{equation} }

It is clear that the four vertices of an isosceles trapezoid lie on a circle
and consequently, its cross ratio has to be real.
To understand the converse we assume that the vertices of a trapezoid lie on
a circle. Then, the four vertices are obtained by intersecting this circle
with two parallel lines. Thus, the trapezoid has a reflection symmetry showing
that it has to be isosceles. This proves the first part of the trapezium lemma.
Since a quadrilateral is uniquely determined by three of its points and its
cross ratio the second part becomes clear by calculating the cross ratio
(\ref{CrossRatioConformal}) of an isosceles trapezoid.

One direction in the proof of the theorem is now obvious:
if the Christoffel transform $F^c$ of an isothermic net $F$ is also a Darboux
transform, then, all the trapezoids spanned by corresponding edges of $F$ and
$F^c$ are isosceles. Thus\footnote{Moreover, after contemplating the relation
between the lengthes of the edges and of the diagonals of an isosceles
trapezoid we see that also the distances $|F_{m,n}-F^c_{m+1,n}|$ and
$|F_{m,n}-F^c_{m+1,n}|$ are constant.\label{MeanCurv1}},
$|F_{m,n}-F^c_{m,n}|=|F_{m+1,n}-F^c_{m+1,n}|=|F_{m,n+1}-F^c_{m,n+1}|.$
The other direction becomes clear, again, by calculating the cross ratio of an
isosceles trapezoid: up to the ``correct'' scaling factor $\lambda^c$ for
the Christoffel transform $F^c$, corresponding edges of $F$ and $F^c$ have
reciprocal lengthes (\ref{ChristoffelTransform}) and consequently, the cross
ratios (\ref{CrossRatioConformal}) become
\begin{equation}\begin{array}{rcr}
	\DV(F_{m,n},F^c_{m,n},F^c_{m+1,n},F_{m+1,n})&=&
		\lambda^c|F_{m,n}-F^c_{m,n}|^2\cr
	\DV(F_{m,n},F^c_{m,n},F^c_{m,n+1},F_{m,n+1})&=&
		-\lambda^c|F_{m,n}-F^c_{m,n}|^2\cr
	\end{array}\end{equation}
--- which proves the theorem.

Now, we are able to give a definition of discrete cmc nets involving the
mean curvature:

\Theorem{Definition (Discrete net of constant mean curvature)}{
An isothermic net $F:\Gamma\to\R^3$ is a {\df net of constant mean curvature
$H$} if there is a Christoffel transform $F^p:\Gamma\to\R^3$ of
$F$ in constant distance\footnote{As before (cf. footnote \ref{SphereCMC}),
the sphere is {\em excluded} by this definition.}
\begin{equation}|F_{n,m}-F^p_{n,m}|^2=\frac{1}{H^2}.\end{equation}
$F^p$ is called the {\df parallel cmc net} of $F$.}

Note that this definition is symmetric in $F$ and its Christoffel transform $F^p$.
So, by definition, the (correctly scaled) Christoffel transform of a discrete
net of constant mean curvature $H$ is a discrete net of constant mean curvature
$H$, too.
Clearly, rescaling of the Christoffel transform results in a discrete net
of (generally) another constant mean curvature.
Also note that the sign of the (constant) mean curvature $H$ is not determined
in the above definition\footnote{Choosing a sign for the constant mean curvature
$H$ would correspond to the choice of a unit normal field $N:\Gamma\to S^2$.
Note, that also the definition of the mean curvature function scetched in
footnote \ref{MeanCurv} includes no sign choice.\par
From our previous discussions (cf. footnote \ref{MeanCurv1}), it is clear that
this mean curvature function (in the sense of footnote \ref{MeanCurv}), indeed,
equals the constant $H$ for a ``discrete net of constant mean curvature $H$''
(in the sense of the above definition).}.

From the work of Bianchi \cite{Bianchi2} it is well known that (smooth) cmc
surfaces allow $\infty^3$ Darboux transforms in cmc surfaces. In \cite{Q3}
we showed that a constant mean curvature Darboux transform $\hat{F}:U\to\R^3$
of a constant mean curvature surface $F:U\subset\R^2\to\R^3$ has constant distance
\begin{equation}\label{CMC}|\hat{F}-F^p|^2=\frac{1-H^cH}{H^2}
	\end{equation}
to its parallel constant mean curvature surface $F^p$.
Herein, the mean curvature $H^c$ of $F^c$ of $F$ is used to determine the
scaling of the Christoffel transform in the Riccati equation (\ref{SmoothRiccati}).
The 3-parameter family mentioned by Bianchi is obtained by solving the Riccati
equation with an initial value satisfying (\ref{CMC}).
A similar theorem holds in the discrete case --- we will obtain it as a
consequence of the following

\Theorem{Supplement to the Hexahedron lemma}{If, in the hexahedron lemma,
two of the faces (quadrilaterals) of the constructed hexahedron are (isosceles)
trapezoids, then their opposite faces are (isosceles) trapezoids, too.}

To prove this supplement let --- without loss of generality --- the initial
quadrilateral $(X_1,X_2,X_3,X_4)$ and the first constructed quadrilateral
$(X_1,X_2,Z_2,Z_1)$ be isosceles trapezoids\footnote{If two opposite
quadrilaterals are trapezoids then there is nothing to prove.}.
%
%
Calculating the cross ratios of these two trapezoids yields
\begin{equation}\begin{array}{rcl}
	\mu&=&\pm\frac{|X_1-X_2|^2}{|X_1-X_4||X_2-X_3|}\cr
	\lambda\mu&=&\pm\frac{|X_1-X_2|^2}{|X_1-Z_1||X_2-Z_2|}\cr
	\end{array}\end{equation}
--- where all sign combinations can occur, depending on whether the
two trapezoids are embedded or not.
The remaining two points $Z_3$ and $Z_4$ have to satisfy the condition
\begin{equation}\DV(X_1,Z_1,Z_4,X_4)=\DV(X_2,Z_2,Z_3,X_3)=\lambda.
	\end{equation}
\begin{figure}\centerline{\input{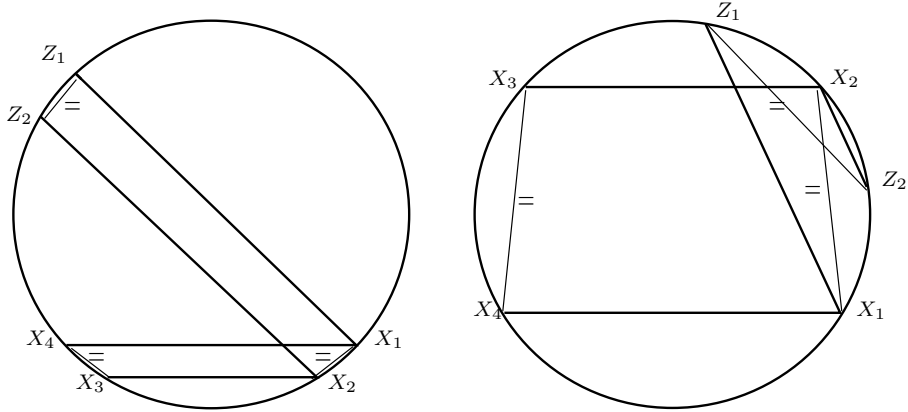}}
	\caption{A supplement to the hexahedron lemma}
	\label{supplement}\end{figure}
In our situation, these two points can be constructed quite explicitely ---
showing that the quadrilaterals $(X_1,Z_1,Z_4,X_4)$ and $(X_2,Z_2,Z_3,X_3)$
have parallel edges:
without loss of generality we may assume that the circles determined by the
point triples $\{X_4,X_1,Z_1\}$ and $\{X_3,X_2,Z_2\}$ have the same radius.
Since $(X_4-X_1)\parallel(X_3-X_2)$ and $(X_1-Z_1)\parallel(X_2-Z_2)$ we
conclude $|X_4-Z_1|=|X_3-Z_2|$ --- this argument is scetched in figure
\ref{supplement}, in case the first two trapezoids are embedded and in
case the initial one is embedded and the second one is not.
Thus, we find a reflection which interchanges $X_4$ with $Z_2$ and $X_3$
with $Z_1$. Now, the two points we had to construct are the images $Z_4$ of
$X_2$ and $Z_3$ of $X_1$ --- it is easily checked that the quadrilaterals
$(X_1,Z_1,Z_4,X_4)$ and $(X_2,Z_2,Z_3,X_3)$ have the correct cross ratio.
Since these quadrilaterals obviously have parallel edges all quadrilaterals
$(X_1,X_2,X_3,X_4)$, $(X_1,X_2,Z_2,Z_1)$, $(Z_1,Z_2,Z_3,Z_4)$ and $(X_4,X_3,Z_3,Z_4)$
are trapezoids and since their cross ratios are real they are isosceles.
This proves the above supplement to the hexahedron lemma.

\IvanOne 

With this knowledge we can now prove the announced

\Theorem{Theorem}{If $F:\Gamma\to\Im\H$ is a discrete net of constant
mean curvature $H$, then any solution $\hat{F}:\Gamma\to\Im\H$ of the
Riccati type system (\ref{Riccati}) with the initial condition
\begin{equation}\label{InitialCMC}
	|\hat{F}_{m_0,n_0}-F^p_{m_0,n_0}|^2=
		\frac{1}{H^2}\left(1-\frac{\lambda}{\lambda^p}\right)
	\end{equation}
is a discrete net of constant mean curvature $H$. Herein, $\lambda^p$ is
the parameter of the Darboux transform which transforms the cmc net $F$
into its parallel cmc net $F^p$.}

To prove this theorem we will have to construct the parallel cmc net
$\hat{F}^p$ of $\hat{F}$. This net is a $\lambda^p$-Darboux transform
of $\hat{F}$ and a $\lambda$-Darboux transform of the parallel cmc net
$F^p$ of $F$ at the same time --- consequently, the construction will
be similar to that in the proof of the permutability theorem for the
Darboux transform (cf. Fig. \ref{BianchiTheorem}).

By the trapezium lemma, the quadrilateral
$(F_{m_0,n_0},F^p_{m_0,n_0},\hat{F}^p_{m_0,n_0},\hat{F}_{m_0,n_0})$
which is first constructed is an isosceles trapezoid because of
(\ref{InitialCMC}).
Since all the quadrilaterals spanned by corresponding edges of $F$ and
its parallel cmc net $F^p$ are also isosceles trapezoids, according to
our supplement to the hexahedron lemma, there is a net
$\hat{F}^p:\Gamma\to\Im\H$ whose edges are parallel to the corresponding
edges of $\hat{F}$ and which has constant distance $\frac{1}{|H|}$ to $\hat{F}$.
The only thing left is to show that $\hat{F}^p$ is a Christoffel transform
of $\hat{F}$ --- which becomes clear by contemplating the cross ratios
\begin{equation}\begin{array}{rcr}
	\DV(\hat{F}_{m,n},\hat{F}^p_{m,n},\hat{F}^p_{m+1,n},\hat{F}_{m+1,n})&=&
		\lambda^p\ \cr
	\DV(\hat{F}_{m,n},\hat{F}^p_{m,n},\hat{F}^p_{m,n+1},\hat{F}_{m,n+1})&=&
		-\lambda^p:\cr
	\end{array}\end{equation}
since $|\hat{F}_{m,n}-\hat{F}^p_{m,n}|^2=\frac{1}{H^2}$ it follows
\begin{equation}\begin{array}{rcr}
	(\hat{F}_{m+1,n}-\hat{F}_{m,n})(\hat{F}^p_{m+1,n}-\hat{F}^p_{m,n})&=&
		\frac{1}{\lambda^pH^2}\,\cr
	(\hat{F}_{m,n+1}-\hat{F}_{m,n})(\hat{F}^p_{m,n+1}-\hat{F}^p_{m,n})&=&
		-\frac{1}{\lambda^pH^2}.\cr
	\end{array}\end{equation}
Thus, $\hat{F}^p$ is a Christoffel transform of $\hat{F}$ and consequently
it is its parallel cmc net in distance $\frac{1}{|H|}$.
This completes the proof.

Now, let's assume we have two cmc Darboux transforms $\hat{F}_1$ and
$\hat{F}_2$ of a discrete net $F$ of constant mean curvature $H$,
and their parallel cmc nets $F^p$, $\hat{F}_1^p$ and $\hat{F}_2^p$
which are $\lambda^p$-Darboux transforms of the original nets.
By the ``fancy'' version of Bianchi's permutability theorem,
the picture gets completed with two discrete isothermic nets $\hat{F}$
and $\hat{F}^p$ --- $\hat{F}$ being a $\lambda_2$-Darboux transform
of $\hat{F}_1$ and a $\lambda_1$-Darboux transform of $\hat{F}_2$ and
$\hat{F}^p$ being a $\lambda^p$-Darboux transform of $\hat{F}$, a
$\lambda_1$-Darboux transform of $\hat{F}_2^p$ and a $\lambda_2$-Darboux
transform of $\hat{F}_1^p$.
Since, moreover, $F^p$ and $\hat{F}_i^p$ are the parallel cmc nets of
$F$ and $\hat{F}_i$ the quadrilaterals $(F,F^p,\hat{F}_1^p,\hat{F}_1)$
and $(F,F^p,\hat{F}_2^p,\hat{F}_2)$ are isosceles trapezoids.
Thus, by our supplement to the hexahedron lemma, the quadrilaterals
$(\hat{F}_1,\hat{F}_1^p,\hat{F}^p,\hat{F})$ and
$(\hat{F}_2,\hat{F}_2^p,\hat{F}^p,\hat{F})$ are isosceles trapezoids, too.
Consequently, the net $\hat{F}^p$ is the parallel cmc net of the cmc net
$\hat{F}$ --- we just proved the following permutability theorem for
cmc Darboux transforms of discrete cmc nets (cf. \cite{Bianchi2}):

\Theorem{Theorem}{If $\hat{F}_1,\hat{F}_2:\Gamma\to\Im\H$ are two cmc
Darboux transforms of a discrete net $F:\Gamma\to\Im\H$ of constant mean
curvature $H$ with parameters $\lambda_1$ and $\lambda_2$, respectively,
then there exists a net $\hat{F}:\Gamma\to\Im\H$ of constant mean curvature
$H$ which is a $\lambda_2$-Darboux transform of $\hat{F}_1$ and a
$\lambda_1$-Darboux transform of $\hat{F}_2$ at the same time.
The nets $F$, $\hat{F}_1$, $\hat{F}_2$ and $\hat{F}$ have constant cross
ratio $\frac{\lambda_2}{\lambda_1}\equiv\DV(F,\hat{F}_2,\hat{F},\hat{F}_1)$.}

\IvanTwo 

\Theorem{Acknowledgment}{We would like to thank Alexander Bobenko for many
fruitful discussions.}

\end{document}